\documentclass {proc-l}

\newtheorem{theo}{Theorem}

\begin {document}

\title [Hydrodynamic approach to constructing solutions of NLS]
{Hydrodynamic approach to constructing solutions of
nonlinear Schr\"odinger equation in the critical case}

\author {O.S. Rozanova}

\address {Department of Differential Equations,
Mathematics and Mechanics Faculty,
Moscow State University,
GSP-2
Vorobiovy Gory
Moscow   119992,
Russia}

\email {rozanova@mech.math.msu.su}

\subjclass {Primary 35Q55; Secondary 35K55}

\date {September 16, 2003}

\keywords {nonlinear Schr\"odinger equation,
hydrodynamic approach, integral functionals, exact solutions,
blow up solutions}

\begin {abstract}
{Proceeding from the hydrodynamic approach, we construct
exact solutions to nonlinear  Schr\"odinger equation with
special properties. The solutions describe collapse, in finite time,
and scattering, over infinite time, of wave packets. They generalize
known blow-up solutions based on the "ground state."}
\end {abstract}

\maketitle

%\documentclass{article}
%\documentclass{amsart}
%\usepackage{amssymb}

%\textheight= 220mm \textwidth=140mm

%\hoffset=-1cm
%\begin {document}

\section{Preliminaries}
Consider the initial value problem for
nonlinear  Schr\"odinger equation (NLS) in  ${\mathbb R}^n$:
$$i\Psi'_t+\Delta \Psi +|\Psi|^\sigma \Psi=0,\qquad
\Psi({\bf x},t):{\mathbb R}^n_x\times {\mathbb R}^+_t\to {\mathbb
C},\eqno(1)$$ $$\Psi(0,{\bf x})=\Psi_0({\bf x})\in H^1({\mathbb
R}^n).\eqno(2)$$

It is well known that the Cauchy problem (1),(2) has locally in
time a solution of class $C([0,T);H^1({\mathbb R}^n)),
\,T\le\infty$ \cite{GinibreVelo1}, \cite{GinibreVelo2},
\cite{Kato}. The Cauchy problem is also locally well-posed in
$L^2({\mathbb R}^n)$\cite{Gasenave8}, and this space is optimal
\cite{Birnir2}.

Moreover, $T=\infty$ for $\sigma<\frac{4}{n},$ where dispersion
dominates, and for $\sigma\ge\frac{4}{n}$ the solution may "blow up"
in  finite time under certain initial conditions, e.g.\cite
{Weinstein},\cite{Glassey}. More exactly, there exist
initial data $\Psi_0({\bf x})$ and a positive constant
$T(\Psi_0({\bf x}))<\infty$, such that
$\displaystyle\lim\limits_{t\to T(\Psi_0({\bf
x}))}\int\limits_{{\mathbb R}^n}|{\bf\nabla}\Psi({\bf
x},t)|^2\,d{\bf x}=\infty.$

The blow-up corresponds to self-trapping of beams in
the laser propagation.

A very good review with references can be found in
\cite{MerleICM}

For solutions to (1),(2) the following quantities are conserved:
$$N[\Psi]=\|\Psi({\bf x},t)\|_{L_2({\mathbb
R}^n)}=\|\Psi_0({\bf x})\|_{L^2({\mathbb R}^n)} =N\quad \mbox{(the
probability)},\eqno(I)$$
$$P[\Psi]=Im \,\int\limits_{{\mathbb R}^n}\Psi({\bf
 x},t){\bf\nabla}\bar\Psi({\bf x},t)\, d{\bf x}=P\quad \mbox{(the
linear momentum)},\eqno(II)$$ and
$$H[\Psi]=\|{\bf \nabla}\Psi({\bf
 x},t)\|_{L^2({\mathbb R}^n)}^2-\frac{2}{\sigma+2}\|\Psi({\bf
x},t)\|_{L^{\sigma+2}({\mathbb R}^n)}^{\sigma+2}=H \quad \mbox{(the
energy)}.\eqno(III)$$

 Let $\sigma=\frac{4}{n}$.  Assume additionally a stronger decay of initial
data at infinity, namely, $|{\bf x}|\Psi_0({\bf x})\in L_2({\mathbb
 R}^n).$ Then the following important identity holds:
 $$M''(t)=8H,\eqno(3)$$ where $M(t)=\int\limits_{{\mathbb
R}^n}|\Psi({\bf x},t)|^2 |{\bf x}|^2 d{\bf x}.$ This identity appears
in many papers, but \cite{Talanov} seems to be the earliest one. It
implies a very simple sufficient condition for the blow-up.
It is easy to see that $M(t)>0$ for the solution of class
$C([0,T);H^1({\mathbb R}^n)),$
\,$T\le\infty$.
At the same time
(3) implies that
$$M(t)=4Ht^2+M'(0)t+M(0).\eqno(4)$$
The blow-up corresponds to vanishing of $M(t).$
It signifies the concentration of the solution support in a set of
zero measure. From elementary algebra arguments we get
that $M(t)\to 0$ at a finite moment of time $T$ if $K\le 0,$ where
$$K=16HM(0)-(M'(0))^2.\eqno(5)$$ This time $T$ is positive
if $H<0$
or $H\ge 0,\,M'(0)<0.$

In \cite{Ogawa} it was proved that for radial initial conditions
$\Psi_0({\bf x})$ with $H<0$ the solution blows up in finite time
in $H^1({\mathbb R}^n)$ with no integrability conditions on
$|\Psi|^2|{\bf x}|^2$ (see also \cite{Nava} in the context).

\section{Hydrodynamic interpretation}
Below we use the hydrodynamic approach due to Madelung, e.g.
\cite{Madelung}. Namely, we represent the solution in
trigonometric form, that is
$\Psi(t,{\bf x})=A(t,{\bf x}) \exp(i\phi(t,{\bf x}))$, where
$A(t,{\bf x})\ge 0$ and $\phi(t,{\bf x})$ are real functions, the
amplitude and the phase of the wave, respectively.

Note that if $A(t,{\bf x})$ is compactly supported, then $\Psi(t,{\bf
x})$ is compactly supported, too. The support does not depend on the
phase function $\phi$.

Substituting this representation in (1) and taking the
real and imaginary parts of the resulting equation, we obtain the
following system
$$A'_t+2(\nabla A,\nabla\phi)+ A\Delta\phi=0,\eqno(6)$$
$$A\phi'_t+A|\nabla\phi|^2={\Delta A}+A^{\sigma+1}.\eqno(7)$$

Further, we multiply (6) by $2A$ and apply the gradient operator
to (7). Denote by $\rho$ the probability density, $A^2=|\Psi|^2$, and
by $V$ double the gradient of the phase function, $2\nabla\phi$.
Thus, the final hydrodynamic form of (1) consists of two equations

$$\rho'_t+\rho{\rm div} V+(V,\nabla \rho)=0,\eqno(8)$$

$$\rho(V'_t+(V,\nabla)V)=2(A\nabla(\Delta
A+A^{\sigma+1})-\nabla A (\Delta
A+A^{\sigma+1})),\quad A=\rho^{1/2}.\eqno(9)$$

The only difference from the traditional gas dynamics is the
"exotic" pressure on the right-hand side of (9). This type of
pressure changes the character of singularity completely, but
allows to use the same methods as in the gas dynamics.

The data
$$\rho(0,{\bf x})=\rho_0({\bf x}),\quad
V(0,{\bf x})=V_0({\bf x})$$ complete the statement of the Cauchy
problem for (8),(9).

It suffices to demand
that $$\rho_0^{1/2}({\bf x})\in H^1({\mathbb R}^n),\quad
\rho_0^{1/2}({\bf x})V_0({\bf x})\in L_2({\mathbb R}^n),  $$
to ensure that the corresponding initial function
$\Psi_0({\bf x})$ belongs to the class  $ H^1({\mathbb R}^n).$

Note also that $\rho$ vanishes as $|{\bf x}|\to\infty$,
however, the same is not necessary for $V$.

The conservation laws (I--III) in the new terms are as follows:
$$\int\limits_{{\mathbb R}^n}\rho d{\bf x}=N
\quad \mbox{(the mass)},\eqno(I')$$ $$\int\limits_{{\mathbb
R}^n}\rho V d{\bf x}={\bf \tilde P}\quad \mbox{(the linear
 momentum)}\eqno(II')$$ and $$\int\limits_{{\mathbb R}^n}(\frac{\rho
 |V|^2}{4}+ \frac{|\nabla\rho|^2 }{4\rho}
-\frac{2}{\sigma+2}\rho^{\sigma/2+1}) d{\bf x}=H \quad
 \mbox{(the energy)}.\eqno(III')$$

One can see that the "kinetic energy" component in $ H$ is the same
as in gas dynamics (up to the multiplier).

Note that in gas dynamics terms we have
$$M(t)=\int\limits_{{\mathbb
 R}^n}\rho({\bf x},t)|{\bf x}|^2 d{\bf x},$$ by virtue of (8)
$$M'(t)=2 \int\limits_{{\mathbb
 R}^n}({\bf x},V({\bf x},t))\rho({\bf x},t) d{\bf x},$$
here $\bf x$ is the radius-vector of point in the space.

Let us introduce one more functional:
$$Q_\Lambda(t)=\int\limits_{{\mathbb R}^n}({\bf x,\Lambda})
\rho({\bf x},t)
d{\bf
x},$$
where ${\bf\Lambda}$ is a constant vector from ${\mathbb R}^n.$
From (8) and the conservation of linear momentum (II') we have
$Q'_{\Lambda}(t)=({\bf \tilde P,\Lambda}):= {\bf \tilde
P_\Lambda}=const,$ and
$$
Q_{\Lambda}(t)={\bf \tilde
P_\Lambda}t+Q_{\Lambda}(0).\eqno(10)$$
From the H\"older inequality we have also
$$Q_{\bf\Lambda}^2(t)\le|{\bf\Lambda}|^2 N M(t).$$

In the domains where the amplitude $A>0,$
instead of  (9) we
consider the equivalent equation
$$V'_t+(V,\nabla)V=2\nabla(\frac{\Delta
A}{A}+A^{\sigma}).$$

Now we use the following idea, recently applied to construct
solutions to the gas dynamics equations(see\cite{Serre},
\cite {Roz}, \cite{RozPas}, \cite{Roz1}).  Namely, let the "velocity
field" have the form $$V=a(t)\bf x,\eqno(11)$$ where $a(t)$ is a
time-dependent function.  Thus, the phase function can be restored as
$$\phi(t,{\bf x})=a(t)\frac{|{\bf x}|^2}{4}+\gamma(t),$$ with an
unknown function $\gamma(t).$ Then, from the (linear in $\rho$)
equation (8) we find the density as $$\rho(t,{\bf
x})=\exp(-n\int\limits_0^t a(\tau) d\tau)\rho_0\left({\bf
x}\exp(-\int\limits_0^t a(\tau) d\tau)\right),$$ or $$A(t,{\bf
x})=\exp(-\frac{n}{2}\int\limits_0^t a(\tau) d\tau) A_0\left({\bf
x}\exp(-\int\limits_0^t a(\tau) d\tau)\right),$$ with $A_0({\bf
x})=|\Psi_0 ({\bf x})|.$

Note that in the critical case, $\sigma=\frac{4}{n}$,
$$
\frac{\Delta
A(t,{\bf x})}{A(t,{\bf x})}+A^{\sigma}(t,{\bf
x})=\exp(-2\int\limits_0^t a(\tau) d\tau) (\frac{\Delta
A_0({\bf \xi})}{A_0({\bf \xi})}+A_0^{\sigma}({\bf\xi})),\eqno(12)$$
where ${\bf\xi}={\bf x}\exp(-\int\limits_0^t a(\tau)
d\tau).$

Further, we have from (7) and (12)
$$(a'(t)+a^2(t))\frac{|{\bf x}|^2}{4}+\gamma'(t)=\exp(-2\int
\limits_0^t a(\tau) d\tau)\left(\frac{\Delta
A_0({\bf \xi})}{A_0({\bf \xi})}+A_0^{\sigma}({\bf \xi})\right),$$
or
$$(a'(t)+a^2(t))\frac{|{\bf
\xi}|^2}{4}+\gamma'(t)\exp(-2\int \limits_0^t a(\tau)
d\tau)=\exp(-4\int \limits_0^t a(\tau)
d\tau)(\frac{\Delta A_0({\bf \xi})}{A_0({\bf \xi})}+A_0^{\sigma}({\bf
\xi})).$$

The variables $t$ and ${\bf\xi}$ can be separated if
$$\gamma'(t)=\gamma_0\exp(-2\int \limits_0^t a(\tau)),\quad
\gamma_0\in{\mathbb R}^1.\eqno(13)$$
In this case
$$(a'(t)+a^2(t)){|{\bf \xi}|^2}=4\exp(-4\int
\limits_0^t a(\tau) d\tau)\left(\frac{\Delta A_0({\bf \xi})}{A_0({\bf
\xi})}+A_0^{\sigma}({\bf \xi})-\gamma_0\right).\eqno(14)$$

It follows from (14) that
$$a'(t)+a^2(t)=4k\exp(-4\int_0^t a(\tau)d\tau),\eqno(15)$$
$$\Delta A_0({\bf x})+A_0^{\sigma+1}({\bf x})=(k|{\bf
x}|^2+\gamma_0)A_0({\bf x}),\eqno(16)$$
where $k$ is a constant.

Thus, we seek a special solution to (1) in the form
$$\Psi(t,{\bf x})=\exp(-\frac{n}{2}\int\limits_0^t a(\tau)
d\tau)A_0\left[{\bf x}\exp(-\int\limits_0^t a(\tau) d\tau)\right]$$
$$
\exp\left(ia(t)\frac{|{\bf
x}|^2}{4}\right)\exp\left(i\gamma_0\int_0^t\exp(-2\int_0^\tau
a(\tau_1) d\tau_1)d\tau\right)\exp (i\theta),\eqno(17)$$ with
$\theta\in{\mathbb R}^1.$

Because the origin is not a particular point, without loss
in generality we may consider the velocity field $V=a(t)({\bf x
-x}_0)$, where ${\bf x}_0$ is an arbitrary fixed point.
Then the solution takes the form
$$\Psi(t,{\bf x})=\exp(-\frac{n}{2}\int\limits_0^t a(\tau)
d\tau)A_0\left[({\bf x}-{\bf x}_0)\exp(-\int\limits_0^t a(\tau)
d\tau)\right] $$ $$ \exp\left(ia(t)\frac{|{\bf x-x}_0|^2}{4}\right)
\exp\left(i\gamma_0\int_0^t\exp(-2\int_0^\tau a(\tau_1)
d\tau_1)d\tau\right)\exp (i\theta).\eqno(17')$$

The decay properties of the solution as $|{\bf x}|\to \infty$ depend
on $A_0({\bf x}).$ The physical sense requires that the solution
should be of the class $L^2({\mathbb R}^n).$
If we wish to consider solutions from the space
$C([0,T);H^1({\mathbb R}^n)),
\,T\le\infty,$ natural for the existence and uniqueness
to the Cauchy problem (1), (2), we have to choose the initial
data
$$\Psi_0({\bf x})=A_0({\bf x}-{\bf x}_0) \exp\left(ia_0\frac{|{\bf
x-x}_0|^2}{4}\right) \exp (i\theta),$$ where $A_0({\bf x})$ is
a non-negative
solution to (14) belonging to $H^1({\mathbb R}^n)$
and $|{\bf x}|A_0({\bf x})\in L_2({\mathbb R}^n).$
Note that with the function $A_0$ of this class we have also
conservation laws (I - III) for solutions (17) and (17').

\section{Time evolution}

Let us investigate the qualitative behavior of $a(t)$ governed by
(15). It is easy to see that in the case  $k<0$ for any initial
datum $a(0)=a_0$ there exists a moment $T_*$ such that
$a(t)\to-\infty,\, t\to T_*.$ Really, as  $a'(t)<\epsilon<0,$ then
under any $a_0$ there exists a moment $T_1\ge 0$ when
$a(T_1)<0.$ Further, the comparison theorem shows that $a(t)\le
\tilde a(t),$ where $\tilde a(t)$ is a solution to the Cauchy
problem $\tilde a(t)'=-\tilde a^2(t),\,\tilde a(T_1)=a(T_1),\,t\ge
T_1.$ This means that if $k<0,$ the solution $\Psi$ to (1)
of the form (17),
with the amplitude $A_0$ satisfying (16), localizes at the origin
within a finite interval of time, provided $N<\infty.$

We can also analyze the  differential corollary of (15)
$$a''(t)+6a(t)a'(t)+4a^3(t)=0,$$
with the initial data
$a(0)=a_0,\,a'(0)=-a_0^2+4k.$
 It gives us, in particular, that in the case $k>0$ the solution
vanishes at infinity as $O(t^{-1}).$

If $k=0,$ any nontrivial solution to (15) blows up at a finite
time $T=-a_0^{-1}.$ The time is positive if $a_0<0.$
Moreover, if $|{\bf x}|\Psi_0({\bf x})\in L^2({\mathbb R}^n),$
we can find $a(t)$ explicitly for any $k.$
Observe that, if $V=a(t)\bf x$, then
$$M'(t)=2a(t)M(t),\eqno(18)$$
therefore
$a(t)=\frac{M'(t)}{2M(t)}.$
The explicit form of $M(t)$ is known, see (4).
Thus,
$$a(t)=\frac{8Ht+M'(0)}{2(4Ht^2+M'(0)t+M(0))}.\eqno(19)$$
If $K=16HM(0)-(M'(0))^2\le 0$ (see (5)), then $M(t)$ tends to zero
($a(t)$ goes to infinity, respectively) within a finite interval
of time.  Moreover, this interval can be readily calculated. If
$K>0,$ then $a(t)\sim t^{-1},\, t\to\infty.$

Now taking into account (15), (18), (19) we can compute
$$k=\frac{K}{16 M^2(0)}=\frac{16HM(0)-(M'(0))^2}{16
M^2(0).}\eqno(20).$$

There are situations where we can express $a(t)$ through
$Q_{\bf\Lambda}(t).$ Namely, for $V=a(t)\bf x$
we get
${\bf\tilde P_\Lambda}= a(t)Q_{\bf\Lambda}(t).$
Therefore, if there is $\bf\Lambda\in{\mathbb
R}^n$ such that $Q_{\bf\Lambda}(t)\ne 0 \,({\bf\tilde P_\Lambda}\ne
0),$ then taking into account (10) we have
$$a(t)=\frac{1}{t+Q_{\bf\Lambda}(0){\bf\tilde P_\Lambda}^{-1}}=
\frac{1}{t+a_0^{-1}}.$$
From (18) we obtain now $M(t)=M(0)a_0^2(t+a_0^{-1})^2.$
Comparing the result with (4) we can see that
$4H= M(0)a^2(0)$ and $K=k=0.$

Note that if $A_0$ is radial, then  $Q_{\bf\Lambda}(t)= 0$
for any $\bf\Lambda.$

\section{Evolution of wave packets}

Summarizing the above results, we can formulate the following
theorem:

\begin {theo}  Suppose that (16) has at least one
nonnegative solution $A_0({\bf x}).$
Then equation (1) has % the
a special solution
given by the explicit formula (17)((17')),
with the function $a(t)$ governed by equation (15).

If
$
{{\bf x}A_0({\bf x})}\in L^2({\mathbb R}^n),$ then the
formula (17') can be written as
$$
\Psi(t,{\bf x})= \left(\frac{M(0)}{M(t)}\right)^{\frac{n}{4}}
A_0\left[\left(\frac{M(0)}{M(t)}\right)^{\frac{1}{2}}
({\bf x}-{\bf x}_0)\right]
$$
$$
\exp\left(i\frac{M'(t)}{8M(t)}|{\bf x}-{\bf x}_0|^2\right)
\exp\left(i\gamma_0 M(0)\int_0^t M^{-1}(\tau)d\tau\right)
\exp (i\theta),\eqno(21)$$
with the quadratic function $M(t)$ having the explicit form (4).

The behavior of the solution depends
on the sign of the constant $k$ (see (20)).

If $k>0,$ then the
solution decays. Namely,
$$\max\limits_{{\bf x}\in{\mathbb
R}^n}|\Psi(t,{\bf x})|=
A^+\left(\frac{M(0)}{4H}\right)^{\frac{n}{4}}t^{-n/2},\quad
A^+=\max\limits_{{\mathbb R}^n}A_0({\bf x}).$$

If $k\le 0$, then the solution blows up at the point ${\bf
x}_0$ at a finite moment of time $T.$
This time is positive in the following cases:

(i)\quad if $H=0,\,M'(0)<0,\,$ then  $T=\frac{-M(0)}{M'(0)};$

(ii)\quad if $H>0,\,M'(0)<0,\,$ then  $T=\frac{-M'(0)-\sqrt
{|K|}}{2H};\quad$

(iii)\quad if $H<0,\,$ then  $T=\frac{-M'(0)+\sqrt
{|K|}}{2H}.$

Moreover, for $k=0,$
$$\max\limits_{{\bf x}\in{\mathbb R}^n}|\Psi(t,{\bf x})|=A^+
\left(\frac{M(0)}{4H}\right)^{\frac{n}{4}}
(T-t)^{-n/2},$$
for $k<0,$
$$\max\limits_{{\bf x}\in{\mathbb R}^n}|\Psi(t,{\bf x})|=A^+
\left(\frac{M(0)}{\sqrt{|K|}}\right)^{\frac{n}{4}}
(T-t)^{-n/4}.$$
\end{theo}

\vskip1cm

{\sc Remark 1.}
If there exists
$\bf\Lambda\in{\mathbb
R}^n$ such that $Q_{\bf\Lambda}(t)\ne 0 \,({\bf\tilde P_\Lambda}\ne
0),$  then only the situation with $k=0$ may be realized.

{\sc Remark 2.}
For the solution of form (19) we have
$$\|\nabla\Psi({\bf x},t)\|_{L^2({\mathbb R}^n)}^2=
\exp(-2\int a(t)dt)
\|\nabla A_0({\bf x})\|_{L^2({\mathbb R}^n)}^2
+$$
$$
\frac{1}{4} a^2(t)\exp(2\int a(t)dt)
\|{\bf x} A_0({\bf x})\|_{L^2({\mathbb R}^n)}^2=$$$$
\frac{M(0)}{M(t)}\|\nabla A_0({\bf x})\|_{L_2({\mathbb
R}^n)}^2+
\frac{(M'(t))^2}{16M(t)}=O(T-t)^{-\lambda},$$
where $\lambda=1$ (the lower estimate for the blow-up order,
see \cite{Thutsumi1})
for $k<0,$
 and $\lambda=2$ for $k=0.$

\section{Comparison with previous results}

In the theory of NLS for the critical case $\sigma=\frac{4}{n}$
the crucial role is played by the so-called ground state, i.e. the
positive radially-decreasing solution
to the elliptic problem
$$
\Delta u+|u|^\sigma u-u=0, \, u\in H^1({\mathbb R}^n).\eqno(22)$$ It
is known \cite{Lions}, that the solution with such properties is
unique and exists at least for $n=1, 2, 3.$ The solution belongs to
$C^2({\mathbb R}^n),$ and $|D^\alpha u|\le C\exp (-\delta|x|),\,$
where $C,\delta$ are positive constants.  We denote it by $R({\bf
x}).$

It is known that $H[R({\bf x})]=0$\cite {Weinstein2}.

It was proved \cite {Weinstein2} that if $\|\Psi_0({\bf
x})\|_{L^2({\mathbb R}^n)}< \|R({\bf x})|_{L^2({\mathbb R}^n)}$,
then the solution to the problem (1), (2) is global in time. If
the solution blows up, then  $\|\Psi_0({|\bf x})\|_{L_2({\mathbb
R}^n)}> \|R({\bf x})\|_{L^2({\mathbb R}^n)}.$ In the case
$\|\Psi_0({\bf x})\|_{L^2({\mathbb R}^n)}= \|R({\bf
x})\|_{L^2({\mathbb R}^n)}$ the solution either blows up or not.

If it does blow up, it necessarily has the following special form
based on the ground state \cite{Merle1,Merle2}: $$\Psi({\bf
x},t)=\exp(i\theta)\exp \left(i(-\frac{\omega^2}{t-T}+\frac{|{\bf
x}-{\bf x}_0|^2}{4(t-T)})\right)
\left(\frac{\omega}{t-T}\right)^{\frac{n}{2}}R\left(\frac{\omega({\bf
x}-{\bf x}_0)}{t-T}-{\bf x}_1\right), \eqno(23) $$ with certain
${\bf x}_0\in{\mathbb R}^n, {\bf x}_1\in{\mathbb R}^n, \theta\in
{\mathbb R},\omega\in {\mathbb R}_+, T\in {\mathbb R}.$ However,
(23) coincides with (21) for $K=k=0$ and $\gamma_0=1.$ Indeed, in
this case $R(x) $ is a solution to (14),
$M(t)=H(t+\frac{M'(0)}{2H})^2=$$
H(t-\left(-\frac{M(0)}{H}\right)^{\frac{1}{2}})^2,$
$\omega^2=\frac{M(0)}{H},\,$$T=\omega=-\frac{M'(0)}{2H}.$ Here we
choose $M'(0)<0 $ to guarantee the positivity of $T.$ Note that
$H>0$ for solutions of the form (23).

It can be readily demonstrated that
 $$H[\Psi]=
\frac{M(0)}{M(t)}(\|\nabla A_0({\bf x})\|_{L_2({\mathbb
R}^n)}^2-\frac{2}{\sigma+2}\|\nabla A_0({\bf
 x})\|_{L_{\sigma+2}({\mathbb R}^n)}^{\sigma+2}+
\frac{(M'(t))^2}{16M(t)}=$$
$$\frac{M(0)}{M(t)}H[A_0]+
\frac{(M'(t))^2}{16M(t)}.$$
Thus, for $A_0=R,$ we get
$H=\frac{(M'(t))^2}{16M(t)}>0.$

It follows from the above results that
if $A_0({\bf x})$ is a solution to (16) with $k<0$
from $H^1({\mathbb R}^n),$ then
$$\|A_0({\bf x})\|_{L^2({\mathbb R}^n)}>\|R({\bf
x})\|_{L^2({\mathbb R}^n)},$$
because in this case the solution (17) blows up
and it is not of the form (23).

The profile of our solution in the case $k<0$ is
different from that of the solution (23) (corresponding to $k=0$),
and so is the rate of blow-up (see Remark 2 of Section 4.)

Note that in the case $k<0$ the solution to equation (16),
considered in the space ${\mathbb R}^n,$ oscillates as $|x|\to
\infty.$ For linearized equation ($n=1$) we can even get the
explicit solution:
$$ u=|x|^{-1/2}(C_1 W_1(\frac{i\gamma_0}{2\sqrt{2k}},\frac{1}{4},
\frac{1}{2}i\sqrt{2k}x^2)+C_2
W_2(\frac{i\gamma_0}{2\sqrt{2k}},\frac{1}{4},
\frac{1}{2}i\sqrt{2k}x^2)),$$ where $W_1, W_2 $ are the Whittaker
functions, $C_1, C_2$ are constants. This highly oscillating
function does not belong to $L^2({\mathbb R}^n).$ So we cannot
hope that the solution to nonlinear perturbed problem (16) is
positive and  belongs to $L^2({\mathbb R}^n) $($H^1({\mathbb
R}^n)).$

However, we can consider (in higher dimensions too) the solution
to (16) given in  $\Omega\subset {\mathbb R}^n.$ To be exact, now
we deal with the Dirichlet problem for (16) with zero boundary
conditions. Regarding this Dirichlet problem, for example, there
is the following result due to \cite {NiNissbaum}.
\begin{theo}[\cite {NiNissbaum}]
Let $\Omega$ be a bounded, smooth domain in ${\mathbb R}^n, n\ge
2,$ and $g:{\mathbb R}^+\times\bar\Omega\to  {\mathbb R}^+$ a
locally Lipshitzian map. Consider  the elliptical boundary problem
$$\Delta u+u^p+\epsilon g(x,u)=0,\,u|_{\partial \Omega}=0.\eqno(*)$$
If $1<p<\frac{n+1}{n-2}$ $(p>1$ for $n=1,2)$, there exists
$\epsilon_0>0$ such that, for $0\le\epsilon<\epsilon_0,$ (*) has a
solution $u=u_\epsilon,$ which is positive on $\Omega.$
\end{theo}

In our situation $g(x,u)=\frac{-k|x|^2-\gamma_0}{\epsilon}u$
satisfies to the theorem condition for $k<0$ in a ball from
${\mathbb R}^n,$ for $\gamma_0<0$ this condition holds for all
${\mathbb R}^n.$

The solution to the Dirichlet problem is classical, that is it
belongs to $C^2(\Omega)\cap C(\bar\Omega),$ therefore it can be
extended to all over the space ${\mathbb R}^n$ at least as
solution from $L^2({\mathbb R}^n).$

\section{Further generalization}

Let us consider the following velocity field:
$${\bf V}=a(t){\bf x}+b(t){\bf \Lambda},\eqno(24)$$
with a constant vector ${\bf \Lambda}.$ Then the phase function
takes the form
$$\phi(t,{\bf x})=a(t)\frac{|{\bf x}|^2}{4}+\frac{1}{2}b(t)
({\Lambda, x})+\gamma(t).$$

Proceeding in the spirit of Section 2, we obtain
$$A(t,{\bf x})=\exp(-\frac{n}{2}\int\limits_0^t
a(\tau) d\tau) A_0\left({\bf x}\exp(-\int\limits_0^t a(\tau)
d\tau)-{\bf\Lambda}\int_0^t
b(\tau)\exp(-\int_0^{\tau}a(\tau_1)d\tau_1)d\tau\right).$$

Denote ${\bf \xi}={\bf x}\exp(-\int\limits_0^t a(\tau)
d\tau)-{\bf\Lambda}\int_0^t
b(\tau)\exp(-\int_0^{\tau}a(\tau_1)d\tau_1)d\tau.$
From (7) and (12) we get
$$\left[a'(t)+a^2(t)\right]\frac{|{\bf
\xi}|^2}{4}+$$
$$
\left[(a'(t)+a^2(t))\int_0^t
b(\tau)\exp(-\int_0^{\tau}a(\tau_1)d\tau_1)d\tau+
(b'(t)+a(t)b(t))\exp(-\int_0^t a(\tau)d\tau)\right]
\frac{({\bf\Lambda,x})}{2}+$$
$$
\left[(a'(t)+a^2(t))(\int_0^t
b(\tau)\exp(-\int_0^{\tau}a(\tau_1)d\tau_1)d\tau)^2
+\right.
$$
$$
\left.2(b'(t)+a(t)b(t))\exp(-\int_0^t a(\tau)d\tau)
\int_0^t
b(\tau)\exp(-\int_0^{\tau}a(\tau_1)d\tau_1)d\tau
+
\right.$$
$$\left.b^2(t)\exp(-2\int_0^t a(\tau)d\tau)\right]
\frac{|{\bf\Lambda}|^2}{4}$$$$
-
\gamma'(t)\exp(-2\int \limits_0^t a(\tau)
d\tau)=\exp(-4\int \limits_0^t a(\tau)
d\tau)\left(\frac{\Delta A_0({\bf \xi})}{A_0({\bf
\xi})}+A_0^{\sigma}({\bf \xi})\right).$$

There are two possibilities for the separation  of variables.

I. The functions $a(t), \gamma(t)$ and $A_0({\bf x}_1),\,{\bf x}_1={\bf
x}+{\bf \Lambda}$ satisfy  equations (15), (13) and (16),
respectively, $b(t)=a(t).$ So we return to the formula
(17') considered above.

II.  The functions $a(t), b(t), \gamma(t)$ and $A_0({\bf x})$
satisfy the following equations:
$$a'(t)+a^2(t)=0,\eqno(25)$$
$$b'(t)+a(t)b(t)=2k_1\exp(-3\int_0^t a(\tau)d\tau),\eqno(26)$$
$$\gamma'(t)=\gamma_0 \exp(-2\int_0^t a(\tau)d\tau)
-
$$
$$
\left[4k_1 \exp(-2\int_0^t a(\tau)d\tau)+
\int_0^t
b(\tau)\exp(-\int_0^{\tau}a(\tau_1)d\tau_1)d\tau
+ b^2(t)\right]
\frac{|{\bf\Lambda}|^2}{4},\eqno(27)
$$
$$\Delta A_0({\bf x})+A_0^{\sigma+1}({\bf x})=(k_1({\bf\Lambda,
x})+\gamma_0)A_0({\bf x}),\eqno(28)$$
where $\gamma_0$ and $k_1$ are constants.

The functions $a(t), b(t), \gamma(t)$ can be found explicitly.

Note that if $k_1=0, \,\gamma_0=1,\,b(t)=0,$ we obtain again a %the
solution that can be represented by formula (23).

The simplest situation is $a(t)=0,\,b(t)=2k_1 t+b_0,\,\gamma(t)=
-\frac{|{\bf\Lambda}|^2}{12 k_1}(2k_1t+b_0)^3+\gamma_0
t+\gamma_1,\,\gamma_0,\gamma_1\in {\mathbb R}^n.$  The
corresponding
solution has the form
$$ \Psi(t,{\bf x})=A_0[{\bf
x}-{\bf\Lambda}t(k_1 t+b_0)]\exp\left(i(k_1t+\frac{b_0}{2})({\bf
\Lambda, x})+\gamma(t) \right),$$ with the function $A_0({\bf x})$
satisfying (28).  For $k_1=0,\,\gamma_0=1,\,b_0=0$ we get the
solitary wave solution $R(t,{\bf x})\exp(i(t+\gamma_1)).$

In the general case $$a(t)=\frac{1}{t+a_0^{-1}},\quad
b(t)=\frac{b_0 a_0^{-1}+2k_1
a_0^2}{t+a_0^{-1}}-2k_1\frac{1}{(t+a_0^{-1})^2},$$
$$\gamma(t)=\int\limits_0^t\left[\frac{\gamma_0}{(\tau+a_0^{-1})^{-2}}
+\frac{k_1|{\bf\Lambda}|^2}
{\tau+a_0^{-1}}\int\limits_0^\tau\frac{b(\tau_1)}{\tau_1+a_0^{-1}}
d\tau_1-\frac{1}{4}|{\bf\Lambda}|^2 b^2(\tau)\right] d\tau.$$

The corresponding solution has the form
$$\Psi({\bf x},t)=\exp(i\gamma(t))\exp
\left(i(b(t)\frac{({\bf x,\Lambda})}{2}
+\frac{|{\bf
x}|^2}{4(t+a_0^{-1})})\right)
$$
$$
\left(\frac{1}{|t+a_0^{-1}|}\right)^{\frac{n}{2}}A_0\left(\frac{{\bf
x}}{t+a_0^{-1}}-{\bf\Lambda}\int\limits_0^t
\frac{b(\tau)}{t+a_0^{-1}}d\tau\right), \eqno(29) $$ with the
function $A_0({\bf x})$ satisfying (28).  The solution (29) has an
interesting feature in the case when $k_1\ne 0,\,A_0({\bf x})\in
L_2({\mathbb R}^n)$  and $a_0<0.$ At a finite time, $T=-a_0^{-1}$,
the mass concentrates at a point which escapes to infinity.

If we suppose that $A_0({\bf x}){\bf x}\in L_2({\mathbb R}^n),$
we can express the solution through the functionals $M(t),
Q_\Lambda (t) $ and $\bf\tilde P.$ Indeed, for the velocity field
(24) we have
$$M'(t)=2a(t)M(t)+b(t)Q_\Lambda(t),\eqno(30)$$
$$\tilde P_\Lambda= a(t)Q_\Lambda(t)
+b(t)N|{\bf\Lambda}|^2.\eqno(31)$$
Further, from (3), (30), (31) we obtain
$$a'(t)=a(t)\frac{2\tilde P_\Lambda Q_\Lambda(t)
-M'(t)N|{\bf\Lambda}|^2}{M(t)N|{\bf\Lambda}|^2-Q^2_\Lambda(t)}-
\frac{\tilde P_\Lambda^2
-8H N|{\bf\Lambda}|^2}{M(t)N|{\bf\Lambda}|^2-Q^2_\Lambda(t)}.$$
Then we use (4) and (10) to get
$$a'(t)=-a(t)\frac{F'(t)}{F(t)}+\frac{C}{F(t)},\eqno(32)$$
where $F(t)=At^2+Bt+C,\,$$A=8H N|{\bf\Lambda}|^2-\tilde P_\Lambda^2,$
$B=M'(0)H N|{\bf\Lambda}|^2-2\tilde P_\Lambda Q_\Lambda(0),$
$C=M(0)N|{\bf\Lambda}|^2-Q_\Lambda^2(0).$
Note that by virtue of the H\"older inequality $C\ge 0.$

Further, we find from (32) and (31)
$$a(t)=\frac{Ct}{F(t)}+a_0,\quad b(t)=\frac{{\tilde
P_\Lambda}(A-C)t^2+({\tilde
P_\Lambda}(B-a_0)-CQ_\Lambda(0))t+({\tilde
P_\Lambda}C-a_0Q_\Lambda(0))}{F(t)N|{\bf\Lambda}|^2}.$$

\vskip1cm

{\bf Acknowledgments.} This work was partially supported by
the Russian Foundation for Basic Researches Award No.
03-02-16263.

\vskip1cm

\bibliographystyle {amsplain}

\end{document}